\documentclass[a4paper, 12pt]{article}
\usepackage{enumerate, theorem}
\usepackage{amsmath, amsfonts, amssymb}
\usepackage[height=22.5cm, width=15cm]{geometry}
\usepackage[all]{xy}
\usepackage{mathrsfs, yfonts}

\DeclareMathOperator{\C}{\mathbb{C}}

\newcommand{\A}{\tilde{\mathcal{A}}}

\newcommand{\parag}[1]{\paragraph{\sc{#1.}}}

\newtheorem{thm}{Theorem}[subsection]
\newtheorem{defn}[thm]{Definition}
\newtheorem{cor}[thm]{Corollary}
\newtheorem{prop}[thm]{Proposition}
\newtheorem{lemma}[thm]{Lemma}

\begin{document}

\title{Holomorphic families of \ $[\lambda]-$primitive themes.}

\author{Daniel Barlet\footnote{Barlet Daniel, Institut Elie Cartan UMR 7502 \ Universit\'e de Lorraine, CNRS, INRIA  et  Institut Universitaire de France,
BP 239 - F - 54506 Vandoeuvre-l\`es-Nancy Cedex.France ;\newline  e-mail: daniel.barlet@univ-lorraine.fr}.}

\date{04/07/13, corrections 28/04/13}

\maketitle

\section*{Abstract}

This article is the continuation of  [B. 13-b]  where we show how the isomorphism class of a \ $[\lambda]-$primitive theme with a given Bernstein polynomial may be characterized by a (small) finite number of complex parameters. We construct here a corresponding locally versal holomorphic deformation of \ $ [\lambda]-$primitive themes for each given  Bernstein polynomial. Then we prove the universality of the corresponding ``canonical family'' in many cases. We also give some examples  where no local universal family exists.\\

\parag{AMS Classification} 32S05, 32S25, 32S40, 32 G 34, 34 M 56.

\parag{Key words} Theme, Vanishing period, Bernstein polynomial, filtered Gauss-Manin system, (a,b)-module, Brieskorn module, isomonodromic deformations.

\newpage

\tableofcontents

\section{Introduction and first examples.}

\subsection{Introduction.}

This article is the continuation of [B.13-b] and the reader will find there the motivations to introduce this notion with a definition and many examples. The present work may be seen as a special case of the study of deformations of isomonodromic germs of ``filtered''regular differential equations corresponding to the case where the monodromy has only one Jordan block.\\

The aim of this article is to introduce holomorphic families of \ $[\lambda]-$primitive themes and to prove that, for any given Bernstein polynomial, the canonical family we construct using the results of [B.13-b] is holomorphic and is locally versal near each point. We show that for many fixed Bernstein polynomials this canonical family is moreover (globally) universal. But we also produce some examples of  \ $[\lambda]-$primitive themes for which there do not exist even a local universal family in a neighbourhood.

\subsection{Definitions.}

Let \ $X$ \ be a reduced  complex space. We shall denote \ $\mathcal{O}_X[[b]]$ \ the sheaf of \ $\C-$algebras on \ $X$ \ associated to the presheaf
$$ U \mapsto \mathcal{O}_X(U)[[b]] .$$
It is a sheaf of \ $\mathcal{O}_X-$algebras. For \ $\mathcal{I} \subset (\mathcal{O}_X)^p$ \ a sub-sheaf of \ $\mathcal{O}_X-$modules (resp. $\mathcal{O}_X-$coherent), we shall denote \ $\mathcal{I}[[b]]$ \ the sub-sheaf of \ $\mathcal{O}_X[[b]]-$modules (resp.  \ $\mathcal{O}_X[[b]]-$coherent) of  \ $(\mathcal{O}_X[[b]])^p$ \ which is generated by \ $\mathcal{I}$.\\

\begin{defn}\label{X-(a,b) 1}
Let \ $X$ \ be a reduced complex space. A sheaf of \ $\mathcal{O}_X-(a,b)-$modules \ $\mathbb{E}$ \ on \ $X$ \ is a  locally free sheaf of finite type of \ $\mathcal{O}_X[[b]]-$modules endowed of a sheaf map
$$ a : \mathbb{E} \to \mathbb{E}  $$
which is \ $\mathcal{O}_X-$linear and continuous for the \ $b-$adic topology of \ $\mathbb{E} $ \ and satisfies the commutation relation \ $a.b - b.a = b^2$.\\
A morphism between two sheaves of \ $\mathcal{O}_X-(a,b)-$modules  is a morphism of sheaves of \ $\mathcal{O}_X[[b]]-$modules which commutes with the respective actions of \ $a$.
\end{defn}

\parag{Example} Let \ $\lambda \in ]0,1]$ \ be a  rational number and let \ $N$ \ be a natural integer. Define
$$ \Xi_{X,\lambda}^{(N)} : = \oplus_{j=0}^N \ \mathcal{O}_X[[b]].e_{\lambda,j} $$
and define \ $a$ \ by induction from the relations\footnote{The reader may think that \ $e_{\lambda,j} = s^{\lambda-1}.(Log\,s)^j\big/j! $.}
\begin{equation*}
  a.e_{\lambda,0} = \lambda.b.e_{\lambda,0} \quad {\rm and} \quad
  a.e_{\lambda,j} = \lambda.b.e_{\lambda,j} + b.e_{\lambda,j-1} \quad {\rm for} \ j \geq 1 
\end{equation*}

\parag{Notation} We shall denote, for \ $\Lambda$ \ a finite subset in \ $]0,1] \cap \mathbb{Q}$,
 \ $$\Xi_{X, \Lambda}^{(N)} : = \oplus_{\lambda \in \Lambda} \quad  \Xi^{(N)}_{X,\lambda}.$$
 An holomorphic map of a reduced complex space \ $X$ \ in \ $\Xi_{\Lambda}^{(N)}$ \ will be, by definition,  a global section on \ $X$ \ of the sheaf \ $\Xi_{X, \Lambda}^{(N)}$.\\
 Remark that we have \ $a.\Xi_{X,\Lambda}^{(N)} \subset b.\Xi_{X, \Lambda}^{(N)}$ \ so they are "simple poles" \ $\mathcal{O}_X-(a,b)-$modules.

\bigskip

Let \ $x \in X$. We have an evaluation map \ $\mathcal{O}_X \to \mathcal{O}_X\big/\mathcal{M}_x \simeq \C_x$ \ where \ $\mathcal{M}_x \subset \mathcal{O}_X$ \ is the sub-sheaf of holomorphic germs which vanish at \ $x$. When \ $\mathbb{E}$ \ is a sheaf of \ $\mathcal{O}_X-(a,b)-$modules on \ $X$, we shall have, in an analoguous way, an evaluation map at \ $x$
$$ \mathbb{E} \to E(x) : = \mathbb{E}\big/\mathcal{M}_x[[b]].\mathbb{E}.$$
Then \ $E(x)$ \ is the fiber at \ $x$ \ of the sheaf \ $\mathbb{E}$. We shall consider a sheaf of \ $\mathcal{O}_X-(a,b)-$modules on \ $X$ \ as a family of (a,b)-modules parametrized by \ $X$.

\begin{defn}\label{holom. 0}
Let \ $X$ \ be a reduced complex space. A holomorphic map \\ $\varphi : X \to \Xi_{ \Lambda}^{(N)}$ \ is  {\bf $k-$thematic} when the following condition is satisfied :
\begin{itemize}
\item The \ $\mathcal{O}_X[[b]]-$sub-module \ $\mathbb{E}_{\varphi}$ \ of \ $\Xi_{X, \Lambda}^{(N)}$ \ generated by the \ $a^{\nu}.\varphi, \nu \in \mathbb{N}$ \ is free of rank \ $k$ \ with basis \ $\varphi, a.\varphi, \dots, a^{k-1}.\varphi$.
\end{itemize}
\end{defn}

For each \ $x \in X$ \ we shall denote \ $E(\varphi(x))$ \ the rank \ $k$ \ theme given as 
$$\mathbb{E}_{\varphi}\big/\mathcal{M}_x[[b]].\mathbb{E}_{\varphi} \simeq \A.\varphi(x) \subset \Xi_{\Lambda}^{(N)}.$$


\begin{lemma}\label{holom. 1}
Let \ $X$ \ be a reduced complex space and \ $\varphi : X \to \Xi_{\Lambda}^{(N)}$ \ a \ $k-$thematic holomorphic map ; then the Bernstein polynomial \ $B_{\varphi(x)}$ \ of \ $E(\varphi(x))$ \ is locally constant on \ $X$.
\end{lemma}

\parag{Proof} We may write, by assumption, 
$$ a^k.\varphi = \sum_{j=0}^{k-1} \ S_j.a^j.\varphi $$
where \ $S_1, \dots , S_{k-1}$ \ are global sections on \ $X$ \ of the sheaf \ $\mathcal{O}_X[[b]]$. As for each \ $x \in X$ \ the theme \ $\A.\varphi(x)$ \ has rank \ $k$, its Bernstein element (see [B.09]) is given by
$$ a^k - \sum_{j=0}^{k-1} \ \sigma_j(x).b^{k-j}.a^j  $$
where \ $\sigma_j(x)$ \ is the coefficient of \ $b^{k-j}$ \ in \ $S_j(x)$. Remark that, when it is not zero,  \ $\sigma_j(x).b^{k-j}$ \ is the initial form of \ $S_j(x)$. But the holomorphic function \ $x \to \sigma_j(x)$ \ takes only rational values\footnote{The Bernstein polynomial \ $B$ \ of \ $E(x)$ \ which is related to its Bernstein element  \ $P \in \A$ \ by the relation \ $B(-b^{-1}.a) = (-b)^{-k}.P$, has its roots in \ $-\Lambda + \mathbb{Z} \subset \mathbb{Q}$. So \ $\sigma_{j(x)}$ \ is rational.}, so it is locally constant.  $\hfill \blacksquare$\\

In the \ $[\lambda]-$primitive case, the sequence \ $\lambda_j + j$ \ is non decreasing, so the fundamental invariants \ $\lambda_1, \dots, \lambda_k$ \  which are determined by the Bernstein element of \ $\A.\varphi(x)$ \ are locally constant.

\parag{Remark} Given a holomorphic map \ $\varphi : X \to \Xi_{\Lambda}^{(N)}$ \ it is not enough to check that for each \ $x \in X$ \ the (a,b)-module \ $E(\varphi(x))$ \ is a  theme of rank \ $k$ \ to have, even locally,  a \ $k-$thematic map, as it is shown by the following example :\\
Let \ $\lambda > 1$ \ a rational number and put for \ $z \in \C$ : 
$$ \varphi(z) : = s^{\lambda-1}.Log\,s + (z + b).s^{\lambda-2} = s^{\lambda-1}.Log\,s + z.s^{\lambda-2} + \frac{1}{\lambda-1}.s^{\lambda-1} .$$
Then the Bernstein element of \ $\A.\varphi(z)$ \  is \ $(a - \lambda.b)(a - \lambda.b)$ \ for \ $z \not= 0$, but for \ $z = 0$ \ the Bernstein element of \ $\A.\varphi(0)$ \ is \ $(a - (\lambda+1).b)(a - \lambda.b)$. We conclude using  the previous lemma, that \ $\varphi$ \ is not \ $2-$thematic.\\

\parag{Example} Let \ $X$ \ be any reduced and irreducible complex space and consider \ 
 $\varphi : X \to \Xi_{\lambda}^{(k-1)}$ \ a holomorphic map such that the coefficient of \ $e_{\lambda,k-1}$ \ is equal to \ $b^n.S$ \ where \ $S$ \ is an invertible element of the algebra \ $\mathcal{O}(X)[[b]]$\footnote{This is equivalent to say that the constant term (constant in \ $b$) of \ $S$ \ is invertible in \ $\mathcal{O}(X)$.}, and such that the \ $b-$valuation of \ $\varphi - b^n.S.e_{\lambda,k-1}$ \ is strictly bigger than \ $n$. Then the sheaf \ $\mathbb{E}_{\varphi} : = \sum_{j=0}^{k-1} \ \mathcal{O}_X[[b]].a^i.\varphi$ \ is free of rank \ $k$ \ on \ $\mathcal{O}_X[[b]]$ \ and stable by \ $a$ :\\
It does not reduce the generality to assume that \ $S = 1$, and in this case define \  
$\psi : = (a - (\lambda+n).b).\varphi$; it satisfies the same hypothesis that \ $\varphi$ \ with \ $k$ \ replaced by \ $k-1$ \ and \ $n$ \ replaced by \ $n+1$. An easy induction allows to conclude.\\
Note that in this kind of example we have \ $p_1 = \dots = p_{k-1} =0 $. So this method constructs only very specific \ $[\lambda]-$primitive themes.

The reader will find a general and systematic  method to produce holomorphic \ $k-$thematic maps in section 1.3 (see corollary \ref{exist. k-thm.}).

\begin{defn}\label{holom. 2}
Let \ $X$ \ be a reduced complex space and let \ $\mathbb{E}$ \ be a sheaf of \ $\mathcal{O}_X-(a,b)-$modules on \ $X$. We shall say that \ $\mathbb{E}$ \ is an holomorphic family of rank \ $k$ \  themes parametrized by \ $X$ \ when the following condition is satisfied :
\begin{itemize}
\item There exists an open covering \ $(\mathcal{U}_{\alpha})_{\alpha \in A}$ \ of \ $X$ \ and for each \ $\alpha \in A$ \ a finite set \ $\Lambda_{\alpha} \subset ]0,1] \cap \mathbb{Q}$ \  a \ $k-$thematic holomorphic map
$$ \varphi_{\alpha} : \mathcal{U}_{\alpha} \to \Xi_{\Lambda_{\alpha}}^{(k-1)} $$
with an isomorphism of  sheaves of \ $\mathcal{O}_X-(a,b)-$modules \ $\mathbb{E}_{\vert \mathcal{U}_{\alpha}} \simeq \mathbb{E}_{\varphi_{\alpha}}$ \ on \ $\mathcal{U}_{\alpha}$.
\end{itemize}
\end{defn}

\parag{Pull back} Of course when we have an holomorphic map \ $f : Y \to X$ \ between two reduced complex space and \ $\mathbb{E}$ \ an holomorphic family of \ $\Lambda-$primitive themes parametrized by \ $X$ \ we may take the pull back of this family by setting for \ $y \in Y$ \ 
 $E(y) : = E(f(y))$.  It is easy to see that the pull back of an holomorphic family is again an holomorphic family, the corresponding sheaf of \ $\mathcal{O}_Y-(a,b)-$modules on \ $Y$ \ being the "analytic" pull back \ $f^*(\mathbb{E}) : = \mathcal{O}_Y \otimes f^{-1}(\mathbb{E})$,  simply because when \ $\varphi : X \to \Xi_{\Lambda}^{(N)}$ \ is an holomorphic and  \ $k-$thematic map, the composition \ $\varphi\circ f : Y \to \Xi_{\Lambda}^{(N)}$ \ is again holomorphic and \ $k-$thematic.

\subsection{Existence of holomorphic k-thematic maps.}

Our purpose here is to show the following proposition.

\begin{prop}\label{existence}
Let \ $X$ \ and \ $T$ \ two complex manifolds and let \   $F : X \to D\times T$ \  be a proper holomorphic map  which is a submersion over \ $D^*\times T$, where \ $D$ \ is a disc with center \ $0$ \ in \ $\C$ \ and \ $D^{*} = D \setminus \{0\}$. Then to the following data 
\begin{enumerate}[i)]
\item A smooth \ $T-relative$  \ $(p+1)-$differential form \ $\omega$ \ on \ $X$ \ such that 
 $$d_{/T}\omega = 0 = d_{/T}F\wedge \omega .$$

  \item  A vanishing $p-$cycle \ $\gamma$ \ in the generic fiber of  the restriction of \ $F$ \ over \ $D \times \{t_0\}$.
\end{enumerate}
we  associate, locally around the generic point of \ $T$, a \ $k-$thematic holomorphic map which defines the family of themes associated to the family parametrized by \ $T$  \ of vanishing periods associated to \ $\omega$ \ and \ $\gamma$.
\end{prop}

Recall that for \ $t_{0}\in T$ \ fixed, the theme defined by \ $\omega_{t_{0}}$ \ and \ $\gamma$ \ is generated by the asymptotic expansion at \ $s = 0$,  with \ $t_{0}$ \ fixed, of the vanishing period  integral 
 $$\int_{\gamma_{s,t_{0}}} \ \omega_{t_{0}}\big/d_{/T}F.$$
As the map \ $F$ \ is a proper submersion over \ $D^{*}\times T$ \ it is a local \ $\mathscr{C}^{\infty}$ \ fibration and this implies that, for any given \ $\gamma$,  we have a well defined  (multivalued in \ $s$) horizontal family \ $\gamma_{s,t}$ \ of \ $p-$cycles in the fibers of \ $F$ \ over \ $D^{*}\times T$.\\
Note that for \ $\gamma$ \ in the spectral sub-space of the monodromy for the eigenvalue \ $exp(2i\pi.\lambda)$ \ the corresponding theme is \ $[\lambda]-$primitive.\\

The proof of this proposition will be given at the end of this section.\\

Recall that \ $\A$ \ is the \ $b-$completion of the algebra \ $\mathcal{A} : = \C< a,b>$ \ with the commutation relation \ $a.b - b.a = b^{2}$ \ (see [B.13-b]), and that an (a,b)-module is a left \ $\A-$module which is free and finite rank on the sub-algebra \ $\C[[b]]$ \ of \ $\A$.\\

We shall begin by an easy  lemma of algebraic geometry over the algebra  \ $Z : = \C[[b]]$.

 \begin{lemma}\label{rg. sem. cont.1}
 Let \ $E$ \ be a regular  (a,b)-module of rank \ $k$. Fix a \ $\C[[b]]-$basis  \ $e_1, \dots, e_k$ \  of \ $E$ \ and consider \ $E$ \ as the affine space  \ $Z^k$ \ over \ $Z$. Then for each natural integer \ $p$ \ the subset  \ $X_p \subset E = Z^k$ \  defined by
 $$ X_p : = \{ x \in E \ / \ rg( \A.x) \leq p \} $$
 is an algebraic subset of \ $E$;  that is to say that there exists finitely many polynomials \ $P_1, \dots, P_N$ \ in  \ $Z[x_1, \dots, x_k]$ \ such that we have
 $$ X_p = \{ x \in Z^k \ / \ P_j(x) = 0 \quad \forall j \in [1,N] \}.$$
 \end{lemma}
 
 \parag{Proof} As \ $E$ \ is regular with rank \ $k$, for each \ $x \in E$ \ the sub-(a,b)-module \ $\A.x$ \ is regular with rank \ $\leq k$. So it is generated as \ $\C[[b]]-$module by \ $x, a.x, \dots a^{k-1}.x$. To write that the rank of \ $\A.x$ \ is \ $\leq p$ \ it enough to write that all  \ $(q,q)$ \ minors of the matrix of these \ $k$ \ vectors in the basis  \ $e_1, \dots, e_k$ \ are zero for all \ $q = p+1$. This gives the polynomials \ $P_1, \dots, P_N$. $\hfill \blacksquare$\\
 
Here is an immediate consequence.
 
 \begin{cor}\label{rg. sem. cont.2}
 Let \ $X$ \ be a reduced complex space and \ $E$ \ be a  regular (a,b)-module of rank \ $k$. Let \ $f : X \to  E$ \ be a holomorphic map\footnote{By fixing a basis \ $e_1, \dots, e_k$ \ of \ $E$ \ on \ $\C[[b]]$ \ this is, by definition,  a global section of the sheaf \ $\mathcal{O}_X[[b]]^k$.}. Then we have a finite  stratification
  $$ X_0 \subset X_1 \subset \dots \subset X_k = X $$
  by closed analytic subsets, such that, for any \ $q \in [1,k]$ \ the subset \ $X_q\setminus X_{q-1}$ \ is exactly the set of \ $x \in X$ \ where the rank of \ $\A.x$ \ is equal to \ $q$.
  \end{cor}
  
  \parag{Remark} The quotient of two holomorphic functions \ $f ,g : D \to \C[[b]]$ \  with \ $g(0) \not= 0$ \ may be well defined for each value of \  \ $z \in D$, and the function \ $f\big/ g$ \ may not be holomorphic on \ $D$. For instance this is the case for \ $z \to \frac{z+b^2}{z+b} $, because a relation like
  $$ z + b^2 = (z+b).(\sum_{j=0}^{\infty} \ a_j(z).b^j) $$
  implies that \  $a_0 \equiv 1$ \ and  \ $a_1 = \frac{-1}{z}$ ! $\hfill \square$
  
   \begin{lemma}\label{mero}
   Let \ $f , g : X \to \C[[b]]$ \ be two holomorphic functions on a irreducible  reduced complex space \ $X$. Assume that \ $g $ \ does not vanish  and that, for each \ $x \in X$ \ the quotient \ $f(x)\big/g(x)$ \ is in \ $\C[[b]]$. Then there exists a Zariski open dense  set \ $X'$ \ in \ $X$ \ such that the map defined by \ $f\big/g$ \ is holomorphic on \ $X'$.
   \end{lemma}
 
 \parag{Proof} We may assume that \ $f \not\equiv 0$ \ on \ $X$. So there exists two  Zariski open dense  sets \ $X_1$ \ and \ $X_2$ \ such that the valuation in \ $b$ \ of \ $f(x)$ \ (resp. of \ $g(x)$) is constant  \ $= k$ \ on \ $X_1$ (resp. is constant \ $= l$ \ on \ $X_2$). Then our assumption implies \ $k \geq l$ \ and it is clear \ that on the Zariski open dense set   \ $X' = X_1 \cap X_2$  \ the map \ $f\big/g$ \ is holomorphic. $\hfill \blacksquare$\\
 
 Remark that the set where \ $g$ \ vanishes is a closed analytic set. So when \ $g$ \ is not identically zero we may again find an open dense set to apply the lemma.\\
 
 \begin{cor}\label{exist. k-thm.}
 Let  \ $f : X \to \Xi^{(N)}_{\Lambda}$ \ be a non identically zero  holomorphic map of an irreducible reduced complex space \ $X$ \ with value in the  (a,b)-module \ $E : = \Xi^{(N)}_{\Lambda}$. Then there exists a dense open set in \ $X$ \ on which the restriction of \ $f$ \ induces  a \ $k-$thematic holomorphic map via \ $x \mapsto \A.x$ \ where \ $k \leq rank(E)$.
 \end{cor}
 
 \parag{Proof} The point is that we may find an open dense set \ $X'$ \ on which the rank of \ $\A.x$ \ is maximal, thanks to the first lemma above. Then we solve a Cramer system with parameter on \ $X'$ \ to find the functions  \ $x \mapsto S_j(x)\in \C[[b]]$ \ which give the relation
  $$ a^k.f(x) = \sum_{j=0}^{k-1} \ S_j(x).a^j.f(x) .$$
  But these functions are "meromorphic". The second lemma above gives then a  open dense set  in \ $X'$ \ on which\ $f$ \  holomorphic and \ $k-$thematic. $\hfill \blacksquare$\\
  
 \parag{Proof of the proposition \ref{existence}} It is  an immediate consequence of the  corollary above that  there exists an open dense set in \ $T$ \ on which the vanishing period associated to \ $\omega$ \ and \ $\gamma$ \ is a holomorphic family of themes, because it is locally \ $k-$thematic.$\hfill \blacksquare$

\parag{Example : the rank 1 case}
Let \ $X$ \ be a connected reduced complex space, and let \ $\varphi : X \to \Xi_{\lambda}^{(0)}$ \ be an holomorphic \ $1-$thematic map. Then there exists \ $S \in \Gamma(X, \mathcal{O}_X[[b]])$ \ which is invertible\footnote{Note that if for some \ $x \in X, S(x)$ \ is not invertible, the Bernstein element jumps, so the map is not \ $1-$thematic.} and an integer\ $n$ \ such that \ $\varphi = S.s^{\lambda+n-1}$. So the sheaf \ $\mathbb{E}_{\varphi}$ \ is isomorphic to the sheaf \ $\mathbb{E}_{\psi}$ \ where \ $\psi : X \to  \Xi_{\lambda}^{(0)}$ \ is the constant map with value \ $s^{\lambda+n-1}$.\\

We study the rank \ $2$ \ case in the next two results.

\begin{prop}\label{holom. rank 2}
Fix \ $\lambda_1 > 1 $ \ a rational number in \ $[\lambda] \in \C\big/\mathbb{Z}$ \ and a positive integer \ $p$. Let \ $X$ \ be a reduced complex space, and let \ $\varphi : X \to \Xi_{\lambda}^{(1)}$ \ be a holomorphic \ $2-$thematic map, such that the fundamental invariants of the themes \ $E(\varphi(x)), x \in X$ \ are given by \ $\lambda_1, p_1 = p \geq 1$. Then there exists an holomorphic map \ $\alpha : X \to \C^*$ \ such that we have
\begin{enumerate}
\item The map \ $\tilde{\psi} : X \to  \Xi_{\lambda}^{(1)}$ \ defined by
\begin{equation*}
\tilde{\psi}(x) : =  \alpha(x).s^{\lambda_1+p-2}.Log\,s + c(\lambda_1,p).s^{\lambda_1-1} \quad {\rm is} \ 2-{\rm thematic}  \tag{@}
 \end{equation*}
 where \ $c(\lambda_1,p) = \frac{-1}{p}.(\lambda_1-1).\lambda_1\dots (\lambda_1+p-2) $.
\item The sheaves of  \ $\mathcal{O}_X-(a,b)-$modules \ $\mathbb{E}_{\tilde{\psi}}$ \ and \ $  \mathbb{E}_{\varphi}$ \ co{\"i}ncide.
 \end{enumerate}
 \end{prop}
 
 The proof will use the following lemma.
 
 \begin{lemma}\label{rank 2}
 Fix \ $\lambda_1 > 1 $ \ a rational number in \ $[\lambda] \in \C\big/\mathbb{Z}$ \ and let \ $p$ \ be  a positive integer. If \ $\rho \in \C^*$ \ and \ $S \in \C[[b]]$ \ are such that
 $$ \chi : = \rho.s^{\lambda_1+p-2}.Log\,s + S.s^{\lambda_1-2} $$
 satisfies \ $(a -(\lambda_1+p-1).b).\chi = (1 + \alpha.b^p).s^{\lambda_1-1}$ \ for some given \ $\alpha \in \C^*$ \ then we have
 $$ \rho =  \alpha/\gamma \quad {\rm and} \quad S = \frac{-1}{p} + z.b^p $$
 for some \ $z \in \C$, where \ $\gamma : = (\lambda_1-1). \lambda_1.(\lambda_1+1) \dots (\lambda_1+p-2)$.
 \end{lemma}
 
 The proof is an exercise left to the reader.\\
 
 \parag{Remark} As \ $\A.\chi$ \ contains \ $\C[[b]].s^{\lambda_1-1}$ \ it contains \ $\chi_z : = \chi -z.b^p.s^{\lambda_1-2}$ \ and both \ $\chi$ \ and \ $\chi_z$ \ have the same annihilator which is the ideal \\
  $$\A.(a - \lambda_1.b).(1 + \alpha.b^p)^{-1}.(a -(\lambda_1+p-1).b).$$
So  we have for each \ $z \in \C$ \ an automorphism of \ $\A.\chi$ \ which sends \ $\chi$ \ to \ $\chi_z$.
 
 \parag{Proof of the proposition} Without lost of generality we may assume that
   \begin{equation*}
 \varphi(x) = s^{\lambda_1 + p -2}.Log\,s + \Sigma(x).s^{\lambda_1-2} \tag{1}
 \end{equation*}
 where \ $\Sigma \in \Gamma(X, \mathcal{O}_X[[b]])$ \  is invertible in \ $\C[[b]]$ \ for each \ $x \in X$. This is consequence of the fact that we must have for each \ $x \in X$ \ a rank \ $2$ \ theme with fundamental invariants \ $\lambda_1, p_1 = p \geq 1$. Then we have
 \begin{align*}
 & (a - (\lambda_1 + p-1).b).\varphi(x) = \frac{s^{\lambda_1+p-1}}{\lambda_1+p-1} + \Sigma(x).s^{\lambda_1-1} + \\
 & \qquad \qquad  + b^2.\Sigma'(x).s^{\lambda_1-2} - (\lambda_1+p-1).b.\Sigma(x).s^{\lambda_1-2} \\
 & \qquad = (b.\Sigma(x)' - p.\Sigma(x) + \gamma.b^{p}).\frac{s^{\lambda_1-1}}{\lambda_1-1}
 \end{align*}
 where \ $\gamma : = (\lambda_1-1)\lambda_1\dots (\lambda_1+p-2) $  \ in order that 
  $$ \frac{s^{\lambda_1+p-1}}{\lambda_1+p-1} = \gamma.b^p.\frac{s^{\lambda_1-1}}{\lambda_1-1} $$ 
 and where \ $\Sigma(x)'$ \ denote the derivative in \ $b$ \ of \ $\Sigma(x) \in \C[[b]]$. So we have \ $(a - (\lambda_1 + p-1).b).\varphi(x) = S(x).s^{\lambda_1-1}$ \ with 
  $$(\lambda_1 - 1).S(x) : = b.\Sigma(x)' - p.\Sigma(x) + \gamma.b^p.$$
So the constant term in \ $b$ \ in \ $S(x)$ \ is equal to \ $(-p/(\lambda_1-1)).\Sigma(x)(0)$ \ and the constant term\ $\Sigma(x)(0)$ \ of \ $\Sigma(x)$ \ is not zero. Put \ $ S(x) : = S_0(x) + S_p(x).b^p + b.\tilde{S}(x) $ \ where \ $S_{0} (x)$ \ is the constant term in \ $b$ \ and where  \ $\tilde{S}(x)$ \ has no term in \ $b^{p-1}$. From our previous remark we see that \ $S_0$ \ is an invertible holomorphic function on \ $X$.\\
 Now let \ $T \in \Gamma(X, \mathcal{O}_X[[b]])$ \ be a solution of the equation
 $$ b.T(x)' - (p-1).T(x) = \tilde{S}(x) \quad \forall x \in X.$$
 Such a \ $T$ \ exists  because \ $\tilde{S}$ \ has no term in \ $b^{p-1}$. Define now \ $\psi : X \to \Xi_{\lambda}^{(1)}$ \ as
 $$ \psi(x) : = \varphi(x) - T(x).s^{\lambda_1-1} .$$
 As \ $\mathbb{E}_{\varphi}$ \ contains \ $\mathcal{O}_X[[b]].s^{\lambda_1-1}$, thanks to the invertibility of \ $S$ \ in \ $\mathcal{O}_X[[b]]$, the equality \ $\mathbb{E}_{\psi} = \mathbb{E}_{\varphi}$ \ is clear. So \ $\psi$ \ is \ $2-$thematic. \\
 Now compute \ $(a - (\lambda_1 + p-1).b).\psi$:
 \begin{align*}
 & (a - (\lambda_1 + p-1).b).\psi(x) = (a - (\lambda_1 + p-1).b).(\varphi(x) - T(x).s^{\lambda_1-1}) \\
 & \qquad  = S(x).s^{\lambda_1-1} - \big[T(x).a + b^2.T(x)'\big].s^{\lambda_1-1} \\
 & \qquad =  (S_0(x) + S_p(x).b^p).s^{\lambda_1-1}
 \end{align*}
 This shows that \ $\tilde{\psi} : = S_0^{-1}.\psi$ \ satisfies the relation
 $$ (a - (\lambda_1 + p-1).b).\tilde{\psi}(x) = (1 + \alpha(x).b^p).s^{\lambda_1-1} $$
 for each \ $x \in X$ \ where \ $\alpha(x) : = S_0^{-1}(x).S_p(x)$. And, of course, \ $\mathbb{E}_{\tilde{\psi}} = \mathbb{E}_{\psi} = \mathbb{E}_{\varphi}$ \ so conditions  2. and 3.  are satisfied. The lemma shows that there exists an holomorphic function \ $z : X \to \C$ \ such that \ $\tilde{\psi}(x) - z(x).b^{p-1}.s^{\lambda_1-1}$ \ has the desired form. But sending \ $\tilde{\psi}$ \ to \ $\tilde{\psi} - z.b^{p-1}.s^{\lambda_1}$ \ is an automorphism of the \ $\mathcal{O}_X-(a,b)-$module \ $\mathbb{E}_{\varphi}$. $\hfill \blacksquare$\\
 
 The case \ $p_1 = p = 0$ \ is given by the following lemma ; as it is a simple variant of the previous proposition, we let the proof to  the reader.
 
 \begin{lemma}\label{rank 2 bis}
 Fix \ $\lambda_1 > 1 $ \ a rational number in \ $[\lambda] \in \C\big/\mathbb{Z}$. Let \ $X$ \ be a reduced complex space, and let \ $\varphi : X \to \Xi_{\lambda}^{(1)}$ \ be an holomorphic \ $2-$thematic map, such that the fundamental invariants of the themes \ $E(\varphi(x)), x \in X$ \ are given by \ $\lambda_1$ \ and \ $ p_1 = 0$. Then if \ $\psi : = (\lambda_1-1). s^{\lambda_1-2}.Log\,s $ \ there is an isomorphism  of \ $\mathcal{O}_X-(a,b)-$modules between the sheaf \ $\mathbb{E}_{\varphi}$ \ and the sheaf \ $\mathbb{E}_{\psi}$ \ corresponding to the constant map with value \ $\psi$.
 \end{lemma}
 
 Note that the sheaf \ $\mathbb{E}_{\psi}$ \ is given as \ $\mathcal{O}_X[[b]].e_1 \oplus \mathcal{O}_X[[b]].e_2$ \ with \ $a$ \ defined by
 $$ a.e_1 = \lambda_1.b.e_1 \quad {\rm and} \quad a.e_2 = (\lambda_1-1).b.e_2 + e_1 .$$
 
 \bigskip
 
 \begin{defn}\label{parametre 1}
 When \ $E$ \ is a rank \ $2$  \ $[\lambda]-$primitive theme with fundamental invariants \ $\lambda_1, p_1 \geq 1$ \ we shall call the {\bf parameter} of \ $E$ \ the  number\ $\alpha \in \C^*$ \ such that \ $E$ \ is isomorphic to \ $\A\big/\A.(a - \lambda_1.b)(1 + \alpha.b^{p_{1}})^{-1}.(a - (\lambda_1+p_1-1).b) $.
 \end{defn}
 
 The following lemma shows that this number is quite easy to detect.
 
 \begin{lemma}\label{parametre 0}
 Let $E$ \ be a rank \ $2$  \ $[\lambda]-$primitive theme with fundamental invariants \ $\lambda_1, p_1 \geq 1$. Let \ $S \in \C[[b]]$ \ such that \ $S(0) = 1$ \ and such that 
  $$E \simeq  \A\big/\A.(a - \lambda_1.b).S^{-1}.(a - (\lambda_1+p_1-1).b).$$
   Then the coefficient of \ $b^{p_1}$ \ in \ $S$ \ is the parameter \ $\alpha$ \ of \ $E$.
 \end{lemma}
 
 \parag{Proof} Let \ $e$ \ be a generator of \ $E$ \ with annihilator the ideal 
  $$\A.(a - \lambda_1.b).S^{-1}.(a - (\lambda_1+p_1-1).b) ,$$
   and put \ $e_2 : = e$ \ and \ $e_1 : = S^{-1}.(a - (\lambda_1+p_1-1).b).e$. Then \ $e_1, e_2$ \ is a \ $\C[[b]]-$base of \ $E$ \ and so we look for a generator \ $\varepsilon$ \  of \ $E$, write
 $$ \varepsilon : = U.e_2 + V.e_1 $$ 
 with \ $U, V \in \C[[b]]$, which is annihilated by \ $(a - \lambda_1.b)(1 + \alpha.b^p)^{-1}.(a - (\lambda_1+p_1-1).b)$. Remark that we know " a priori" that such an \ $\varepsilon$ \ exists thanks to the proposition \ref{holom. rank 2} with \ $X = \{pt\}$. Then compute 
 \begin{align*}
 & (a - (\lambda_1+p_1-1).b).\varepsilon = b^2.U'.e_2 + U.S.e_1 + b^2.V' - (p_1-1).b.V.e_1.
 \end{align*}
 As we are in a theme the kernel of \ $a - \lambda_1.b$ \ is \ $F_1 = \C[[b]].e_1$. So this implies that \ $U'  = 0$ \ and \ $U = U(0)$. We obtain also the equation
 $$ U(0).S + b^2.V' - (p_1-1).b.V = \rho.(1 + \alpha.b^{p_1}) $$
 for some \ $\rho \in \C$. As we assume \ $S(0) = 1$ \ this implies \ $U(0) = \rho$. The fact that \ $\varepsilon$ \ is a generator of \ $E$ \ insure that \ $U(0) = \rho \not= 0$, and then
 $$ b^2.V' - (p_1-1).b.V = \rho(1 - S) + \rho.\alpha.b^{p_{1}} $$
 and the existence of a solution \ $V \in \C[[b]]$ \ for this equation implies that the coefficient of \ $b^{p_1}$ \ in \ $S$ \ is equal to \ $\alpha$. $\hfill \blacksquare$
 
 \parag{Exercise} Let \ $E$ \ be a rank \ $2$  \ $[\lambda]-$primitive theme with fundamental invariants \ $\lambda_1, p_1= p  \geq 1$. Show that, for \ $\delta \in \mathbb{Q}$ \ such that \ $\lambda_1 + \delta >1$ \  the parameter of \ $E \otimes_{a,b} E_{\delta}$ \ is the same than the parameter of \ $E$.\\
 Deduce then that for \ $-\lambda_1 - p +  \delta + 1 > 1$ \ the parameter of the theme \ $E^* \otimes_{a,b} E_{\delta}$ \ is,  up to sign,  the same than the parameter of \ $E$.
 
 \subsection{Holomorphy criterion.}
 
 Let us go back to holomorphic families of \ $[\lambda]-$primitive themes of any rank.
 
 \begin{prop}\label{holom. 2}
 Let \ $X$ \ be a connected reduced complex space and let \ $\mathbb{E}$ \ be an holomorphic family of rank  \ $k$ \ $[\lambda]-$primitive themes parametrized by \ $X$. Denote \ $\lambda_1, p_1, \dots, p_{k-1}$ \ the corresponding  fundamental invariants. For each \ $j \in [0,k]$ \ there exists an unique  holomorphic family \ $\mathbb{F}_j$ \ of rank \ $j$ \ \ $[\lambda]-$primitive themes parametrized by \ $X$ \ with the following properties :
 \begin{enumerate}[i)]
 \item \ $\mathbb{F}_j \subset \mathbb{F}_{j+1}$, $\mathbb{F}_k = \mathbb{E}$ ;
 \item for each \ $x \in X$ \ the theme \ $F_j(x)$ \ is the normal sub-module of rank \ $j$ \ of \ $E(x)$;
 \item the family \ $\mathbb{E}\big/\mathbb{F}_j$ \ is holomorphic for each \ $j \in [0,k]$.
 \end{enumerate}
 \end{prop}
 
 \parag{Proof} The statement is local on \ $X$ \ and we may assume that \ $\mathbb{E} = \mathbb{E}_{\varphi}$ \\
  where \ $\varphi : X \to \Xi_{\lambda}^{(k-1)}$ \ is a \ $k-$thematic holomorphic map. Up to change \ $\varphi$ \ by the action of an invertible section of the sheaf \ $\mathcal{O}_X[[b]]$ \ we may assume that \\
   $\varphi = s^{\lambda_k-1}.(Log\,s)^{k-1} + \theta $ \ where \ $\theta$ \ is a section of the sheaf \ 
    $\Xi_{X, \lambda}^{(k-2)}$.\\
     Put \ $\psi : = (a - \lambda_k.b).\varphi$. Then \ $\psi$ \ is \ $(k-1)-$thematic, because \ $\psi, a.\psi, \dots, a^{k-2}.\psi$ \ is \ $\mathcal{O}(X)[[b]]-$free and generates \ $\mathbb{E}_{\psi}$ :\\
 If we have \ $\sum_{j=0}^{k-2} \ U_j.a^j.\psi = 0 $ \ with \ $U_j \in \mathcal{O}(X)[[b]]$ \ this implies
 $$ \sum_{j=0}^{k-2} \ U_j.a^{j+1}.\varphi -  \lambda_k.\sum_{j=0}^{k-2} \ U_j.a^j.b.\varphi = 0. $$
 As, by hypothesis \ $\varphi, a.\varphi, \dots, a^{k-1}.\varphi$ \ is \ $\mathcal{O}(X)[[b]]-$free, this gives successively \ $U_{k-2} = 0, U_{k-3} = 0, \dots, U_1 = 0$.\\
 Then it is easy to see that \ $\mathbb{F}_{k-1} : = E_{\psi}$ \ induces \ $F_{k-1}(x)$ \ for each \ $x \in X$. This proves i) and ii).\\
To prove iii) it is enough, by an easy induction, to prove it for \ $j = 1$. In this case, it is sufficient to prove that the composition \ $\theta : = f_{\lambda}\circ \varphi : X \to \Xi_{\lambda}^{(k-2)}$, where \ $f_{\lambda} :  \Xi_{\lambda}^{(k-1)} \to  \Xi_{\lambda}^{(k-2)}$ \ is the quotient by \ $\Xi_{\lambda}^{(0)}$, is \ $(k-1)-$thematic. Remark first that  for \ $x \in X$ \ we have 
 $$F_1(x) = Ker\,f_{\lambda} \cap \A.\varphi(x).$$
  Then we want to prove that if we have \ $\sum_{j=0}^{k-2} \ S_j(x).a^j.\varphi(x) \in \Xi_{\lambda}^{(0)}$, then  
  $$S_j(x) = 0 \quad \forall j \in [0,k-2].$$
   If all \ $S_j(x)$ \ are not zero, there exists \ $q \in \mathbb{N}$ \  and \ $T \in \C[[b]]$ \ invertible, such that
 $$ \sum_{j=0}^{k-2} \ S_j(x).a^j.\varphi(x)  = T.s^{\lambda+q-1} .$$
 But then, \ $(a - (\lambda+q).b).T^{-1}.\big(\sum_{j=0}^{k-2} \ S_j(x).a^j\big) $ \ is a polynomial in \ $a$ \ with coefficients in \ $\C[[b]]$, of degree \ $\leq k-1$, which annihilated \ $\varphi(x)$, contradicting the fact that \ $E(x) = \A.\varphi(x)$ \ has rank \ $k$. Now \ $\theta$ \ is \ $(k-1)-$thematic and induces the family \ $E(x)\big/F_1(x), x \in X$. $ \hfill \blacksquare$\\

 \begin{thm}\label{crt. hol.}
 Let \ $X$ \ be a reduced complex space and let \ $E(x)_{x \in X}$ \ be a family of rank \ $k \geq 2$ \ \ $[\lambda]-$primitive themes with fixed fundamental invariants \ $\lambda_1, p_1, \dots, p_{k-1}$ \ defined by a \ $\mathcal{O}_X-(a,b)-$module \ $\mathbb{E}$. Let \ $\alpha : X \to \C^*$ \ the function which associates to \ $x \in X$ \ the parameter of the rank \ $2$ \ $[\lambda]-$primitive theme \ $E(x)\big/F_{k-2}(x)$. Then the family is holomorphic if and only if the following two conditions are satisfied:
 \begin{enumerate}[i)]
 \item The family \ $F_{k-1}(x)$ \ is holomorphic, that is to say that there exists a \ $\mathcal{O}_X-(a,b)-$sub-module \ $\mathbb{F}_{k-1}$ \ of \ $\mathbb{E}$ \ which is an holomorphic family and take the value \ $F_{k-1}(x)$ \ for each \ $x \in X$.
 \item The function \ $\alpha$ \ is holomorphic on \ $X$.
 \end{enumerate}
 \end{thm}
 
 \parag{Remark} Thanks to proposition \ref{holom. rank 2} the condition ii) of the previous theorem is equivalent to the fact that the family \ $E(x)\big/F_{k-2}(x), x \in X$ \ is holomorphic. So this result allows, using induction, to reduce the problem of the holomorphy of a family to the rank \ $2$ \ case.\\
 
 The proof will use the following lemma.
 
 \begin{lemma}\label{holom. 3}
 Let \ $j, q $ \ be natural integers and \ $\lambda$ \ be a rational number in \ $]0,1]$. Denote \ $H(j,q)$ \ the hyperplan of \ $\Xi_{\lambda}^{(j)}$ \ corresponding to the annulation of the coefficient \footnote{We use here the notations of the example given at the begining of section 1.2.}of \ $b^q.e_{\lambda, j}$. Then the \ $\C-$linear map 
 $$f :  H(j,q)\oplus \C.b^q.e_{\lambda, j+1} \longrightarrow b.\Xi_{\lambda}^{(j)} $$
given by applying \ $(a - (\lambda+q).b)$ \ to each term of the sum and adding,  is an isomorphism of  Frechet spaces. So its inverse is linear and continuous.
 \end{lemma}
 
 \parag{Proof} The following equality is easy, for any \ $(h,m) \in \mathbb{N}$ \ with the convention \ $e_{\lambda,-1} = 0$
 $$ (a - (\lambda+q).b).b^m.e_{\lambda,h} = (m-q).b^{m+1}.e_{\lambda, h} + b^{m+1}.e_{\lambda, h-1} .$$
 It implies that \ $(a - (\lambda+q).b).b^q.e_{\lambda,j+1} = b^{q+1}.e_{\lambda,j} $ \ is in \ $ b.\Xi_{\lambda}^{(j)}$. So the image of \ $\Xi_{\lambda}^{(j)}$ \ by \ $ (a - (\lambda+q).b)$ \ is the hyperplane of \ $b.\Xi_{\lambda}^{(j)}$ \  given by the annulation of the coefficient of \ $b^{q+1}.e_{\lambda,j}$ \ and its kernel is \ $\C.b^q.e_{\lambda, 0}$. The conclusion is then easy. $\hfill \blacksquare$
 
 \parag{Remark} If we begin with an element in \ $b.\Xi_{\lambda}^{(j)}$ \ for which the coefficient of \ $b^q.e_{\lambda, j}$ \ is \ $\rho$, then the coefficient of \ $b^q.e_{\lambda, j+1}$ \ in its image by the inverse map will be also \ $\rho$. In particular, it will be non zero for \ $\rho \not= 0$.
 
 \parag{Proof of the theorem \ref{crt. hol.}} The statement is local on \ $X$ \ so we may assume that we have a \ $(k-1)-$thematic holomorphic map \ $\psi : X \to \Xi_{\lambda}^{(k-2)}$ \ such that \ $\mathbb{E}_{\psi}$ \ defines the family \ $F_{k-1}(x), x \in X$. We may also assume that the map
 $$ \psi - b^{\lambda_{k-1} - \lambda}.e_{k-2} $$
 takes its values in \ $\Xi_{\lambda}^{(k-3)}$, up to a change of  $\psi$ \ by multiplication by an invertible element of \ $\mathcal{O}(X)[[b]]$. Put \ $q : = \lambda_k - \lambda, \ S_{k-1} : = 1 + \alpha.b^{p_{k-1}} $ \ and define
 $$ \varphi : X \to \Xi_{\lambda}^{(k-1)} $$
 as the composition of \ $S_{k-1}.\psi$ \ with the inverse map constructed in the previous lemma, with \ $j : = k-2$.  and the obvious inclusion of \ $ H(k-2, \lambda_k - \lambda) \oplus \C.b^q.e_{\lambda, k-1}$ \ in \ $\Xi_{\lambda}^{(k-1)}$. Note that we have \ $\psi(x) \in b.\Xi_{\lambda}^{(k-2)}$ \ for each \ $x \in X$ \ because of the inclusion  \ $F_{k-1}(x) \subset a.E(x) + b.E(x)$, of the fact that \ $\Xi_{\lambda}$ \ has a simple pole  and also because any \ $\A-$linear map from \ $F_{k-1}(x)$ \ to \ $\Xi_{\lambda}$ \ may be extended  to \ $E(x)$ \ (see [B.05]).\\
 Remark that the coefficient of \ $b^q.e_{\lambda, k-2}$ \ in \ $S_{k-1}.\psi$ \ co{\"i}ncides with the coefficient of \ $b^{p_{k-1}}$ \ in \ $S_{k-1}$, so it is given by \ $\alpha$. As it is non zero, the coefficient of \ $ b^q.e_{\lambda,k-1}$ \ in \ $\varphi$ \ does not vanish, thanks to the relation \ $\lambda_k - \lambda + 1 = \lambda_{k-1} - \lambda + p_{k-1}$. This is of course necessary in order that \ $\A.\varphi(x)$ \ will be a rank \ $k$ \ theme. Now it is easy to see that the holomorphic \ $k-$thematic map \ $\varphi$ \ is such that \ $\mathbb{E}_{\varphi}$ \ induces the family \ $E(x), x \in X$. $\hfill \blacksquare$
 
 \subsection{Holomorphy and duality.} 
 
 The last theorem of this paragraph shows that the twisted duality preserves the holomorphy of a family of \ $[\lambda]-$primitive themes.
 
 \begin{thm}\label{holom. dual}
  Let \ $X$ \ be a reduced complex space and let \ $E(x), x \in X$, be a holomorphic family of  \ $[\lambda]-$primitive themes parametrized by the reduced complex space \ $X$. Let \ $\delta \in \mathbb{Q}$ \ such that for each \ $x \in X$ \ the (a,b)-module \ $E(x)^* \otimes_{a,b} E_{\delta}$ \ is a theme. Then the family  \ $E(x)^* \otimes_{a,b} E_{\delta}, \ x \in X$,  is   holomorphic.
 \end{thm}
 
 The proof will use the next lemma, which is an easy exercise left to the reader
 
 \begin{lemma}\label{holom. decal.}
 Let \ $E(x), x \in X$ \ be an holomorphic family of  \ $[\lambda]-$primitive themes parametrized by a reduced complex space. Let \ $\delta \in \mathbb{Q}$ \ such that for each \
  $x \in X$ \ the (a,b)-module \ $E(x) \otimes_{a,b} E_{\delta}$ \ is a theme. Then the family \ $E(x) \otimes_{a,b} E_{\delta}, x \in X$ \ is holomorphic.
 \end{lemma}
 
 \parag{Proof of the theorem \ref{holom. dual}} We make an induction on the rank \ $k$. As the case \ $k = 1$ \ is trivial and the case \ $k = 2$ \ is already known (see the exercise following the definition \ref{parametre 1}), assume that the result is proved for \ $k-1 \geq 2$ \ and consider the situation in rank \ $k$. Let \ $F_1(x), x \in X$ \ be the family of normal sub-themes of rank \ $1$ \ of the themes \ $E(x), x \in X$. The proposition \ref{holom. 2} implies the holomorphy of the family \ $E(x)\big/F_1(x)$. This proposition gives also the holomorphy of the family \ $F_2(x), x \in X$. The induction hypothesis implies the holomorphy of the family \ $(E(x)\big/F_1(x))^*\otimes_{a,b} E_{\delta}$ \ and also of the family \ $(F_2(x))^*\otimes_{a,b} E_{\delta}$. But now we may apply the theorem \ref{crt. hol.} to the family \ $(E(x))^*\otimes_{a,b} E_{\delta}$ \ because the corresponding family of sub-themes of rang \ $k-1$ \ is the holomorphic family \ $(E(x)\big/F_1(x))^*\otimes_{a,b} E_{\delta}$ \  and the corresponding  family of rank \ $2$ \ quotients is the holomorphic family \ $(F_2(x))^*\otimes_{a,b} E_{\delta}$. $\hfill \blacksquare$

 \section{The canonical family.}
 
 \subsection{Definition and holomorphy.}
 
 We shall fix in this paragraph an integer \ $k \geq 1$ \  a rational number \ $\lambda_1 > k-1$ \ and natural integers \ $p_1, \dots, p_{k-1}$. For \ $j \in [1,k-1]$ \ we defined in [B.13-b] the complex vector space \ $V_j \subset \C[[b]]$ \ in the following way :
\begin{enumerate}
\item If \ $p_j+\dots +p_{k-1} < k-j $ \ put \ $V_j : = \oplus_{i=0}^{k-j-1} \C.b^i $.
\item If  \ \ $p_j+\dots +p_{k-1} \geq  k-j $ \ put \ $q_j : = p_j + \dots + p_{j+h}$ \ where \ $h \in \mathbb{N}$ \ is minimal with the property that \ $q_j \geq k-j$ \ and put \ $ V_j : = \oplus_{i=0}^{k-j-1} \C.b^i \oplus \C.b^{q_j}$.
\end{enumerate}

Define now the affine open set \ $W_j \subset V_j$ \ of  an affine hyperplane of the  complex  vector space \ $V_j$ :
 $$ W_j : = \{ S_j \in V_j \ / \  S_j(0) = 1 \quad {\rm and \ the \ coefficient \ of} \quad b^{p_j} \quad {\rm is} \quad \not= 0 \} .$$
 Then put
 \begin{equation*}
 \mathcal{S}(\lambda_1, p_1, \dots, p_{k-1}) : = \{(S_1, \dots, S_{k-1}) \in \C[[b]]^{k-1} \ / \  S_j \in W_j \quad \forall j \in [1,k-1] \}. \tag{*}
 \end{equation*}
 Note that \ $ \mathcal{S}(\lambda_1, p_1, \dots, p_{k-1})$ \ is always  isomorphic to a product \ $\C^{m}\times (\C^{*})^{n}$. So it is a connected complex manifold.
 
 \begin{defn}\label{fam.can. 1}
 For \ $\sigma \in \mathcal{S}(\lambda_1, p_1, \dots, p_{k-1}) $ \ we define \ $E(\sigma)$ \ as the \ $[\lambda]-$primitive theme (with fundamental invariants \ $\lambda_1, p_1, \dots, p_{k-1}$) \ $\A\big/\A.P(\sigma)$ \ where 
 \begin{equation*}
 P(\sigma) : = (a - \lambda_1.b)S_1^{-1} \dots S_{k-1}^{-1}.(a - \lambda_k.b) \tag{**}
 \end{equation*}
  with  \ $\lambda_{j+1} = \lambda_j + p_j - 1$ \ for \ $j \in [1,k-1]$.
 \end{defn}
 
 \parag{Examples} For \ $k = 1$ \ and \ $\lambda_1$ \ given, the canonical family reduces to the theme \ $E_{\lambda_1}$.
 For \ $k = 2$ \ and \ $\lambda_1, p_1$ \ given, we have
 \begin{enumerate}
 \item For \ $p_1 = 0$,  $\mathcal{S}(\lambda_1, 0) = \{1\}$ \  and the corresponding value of \ $P$ \ is
  $$P : = (a - \lambda_1.b).(a - (\lambda_1-1).b).$$
  
 \item  For \ $p_1 \geq 1$, we have \ $\mathcal{S}(\lambda_1,p_1) = \{ 1 + \alpha.b^{p_1}, \alpha \in \C^* \} \simeq \C^*$ \ and the value of \ $P$ \ corresponding to \ $\alpha \in \C^*$ \ is given by
  $$ P(\alpha) = (a - \lambda_1.b)(1 + \alpha.b^{p_1})^{-1}(a - (\lambda_1 + p_1 - 1).b) .$$
  \end{enumerate}
  
  \begin{thm}\label{fam. can. 2}
  For any given integer \ $k \geq 1$, any rational number \ $ \lambda_1 > k-1$ \ and any natural integers \ $p_1, \dots, p_{k-1}$,  the canonical family is holomorphic.
  \end{thm}
  
  \parag{Proof} The proof will use theorem \ref{crt. hol.} and an induction on \ $k$. First remark that  the family \ $E(\sigma)\big/F_1(\sigma)$ \ is the pull back of  the canonical family parametrized by \ $\mathcal{S}(\lambda_2, p_2, \dots, p_{k-1})$ \ via the obvious projection ( recall that \ $\lambda_2 : = \lambda_1 + p_1 - 1> k-2$ \ as \ $\lambda_1 > k-1$)
  $$ \mathcal{S}(\lambda_1, p_1, \dots, p_{k-1}) \to \mathcal{S}(\lambda_2, p_2, \dots, p_{k-1}). $$
 So  it is a holomorphic family, thanks to our induction hypothesis. Moreover the family \ $F_2(\sigma)$ \ is also the pull back  via the obvious projection 
   $$ \mathcal{S}(\lambda_1, p_1, \dots, p_{k-1}) \to \mathcal{S}(\lambda_1, p_1)$$ 
    and the example above \ gives the holomorphy. So the theorem \ref{holom. dual} allows now to use the theorem \ref{crt. hol.} and this  gives the conclusion . $\hfill \blacksquare$
 
 \subsection{Versality and universality.}
 
 Again we shall fix an integer \ $k \geq 1$, a rational number \ $ \lambda_1 > k-1$ \ and natural integers \ $ p_1, \dots , p_{k-1}$ \ in this paragraph.
 
 \begin{defn}\label{versal. univ. 1}
 Let \ $X$ \ be a reduced complex space, \ $x_0$ \ a point in \ $X$, and let \ $\mathbb{E}$ \ an holomorphic family of  \ $[\lambda]-$primitive themes with fundamental invariants \ $\lambda_1, p_1, \dots, p_{k-1}$. We shall say that this family is {\bf versal near \ $x_0$} \  when, for any reduced complex space \ $Y$ \ with a base point \ $y_0$ \  and any holomorphic family  \ $\mathbb{H}$ \ of  \ $[\lambda]-$primitive themes with fundamental invariants \ $\lambda_1, p_1, \dots, p_{k-1}$, parametrized by \ $Y$, such that \ $H(y_0)$ \ is isomorphic to \ $E(x_0)$, there exists an holomorphic map \ $ f : U \to X$ \ of an open neighbourhood \ $U$ \ of \ $y_0$ \ in \ $Y$ \ and an isomorphism of \ $\mathcal{O}_U-(a,b)-$modules \ $\theta : f^*(\mathbb{E}) \simeq  \mathbb{H}_{\vert U}$.\\
 We shall say that \ $\mathbb{E}$ \ is {\bf universal near \ $x_0$}  \  when  the germ at \ $y_0$ \ of  such an holomorphic map \ $f$ \  is always unique.
 The family will be called {\bf  locally versal}  when it is versal  in a neighbourhood of each point. The family will be called {\bf universal} when it is locally universal and when each isomorphy class appears exactly one time in the family.
 \end{defn}
 
 Of course a versal family contains all isomorphy class of  \ $[\lambda]-$primitive themes with fundamental invariants \ $\lambda_1, p_1, \dots, p_{k-1}$. Any  isomorphy class appears exactly one time in an universal family. The existence of an universal family is the same problem that the representability of the functor which associates to a reduced complex space \ $Y$ \ the set of holomorphic families of   \ $[\lambda]-$primitive themes with fundamental invariants \ $\lambda_1, p_1, \dots, p_{k-1}$, parametrized by \ $Y$. \\
 The existence of a versal family implies the finiteness of the number of complex parameters in order to determine an isomorphy class of such a theme.
 
 \begin{thm}\label{versal univ. 2}
 For any choice of an integer \ $k \geq 1$, and of fundamenetal invariants given by  a rational number 
  $ \lambda_1 > k-1$ \  and  natural integers \ $ p_1, \dots , p_{k-1}$ \ the canonical family parametrized by \ $\mathcal{S}(\lambda_1, p_1, \dots, p_{k-1})$ \ is locally  versal.\\
 If any theme corresponding to  our choice of \ $k \geq 1$, $ \lambda_1 > k-1$ \  and \ $ p_1, \dots , p_{k-1}$ \ has an unique canonical form, then this canonical family is (globally) universal.
 This is  the case, for instance, when \ $p_j \geq k-1 \quad \forall j \in [1,k-1]$.
 \end{thm}
 
 \parag{Proof} First we shall show the versality by induction on the rank \ $k \geq 1$. The cases \ $k = 1$ \ and \ $k = 2$ \ are already known (see above section 1.2.), so we shall assume that \ $k \geq 3$ \ and the case \ $k-1$ \ known. \\
 Let \ $\mathbb{E}$ \ be an holomorphic family of  \ $[\lambda]-$primitive themes with fundamental invariants \ $\lambda_1, p_1, \dots, p_{k-1}$ \ parametrized by the reduced complex space \ $X$. Fix \ $x_0$ \ in \ $X$ \ and denote \ $\mathbb{F}_1 \subset \mathbb{E}$ \ the sub-sheaf of \ $\mathcal{O}_X-(a,b)-$modules \ defining the family of normal sub-themes of rank 1 in \ $\mathbb{E}$. Then \ $\mathbb{E}\big/\mathbb{F}_1$ \ is an holomorphic family of  \ $[\lambda]-$primitive themes with fundamental invariants \ $\lambda_2, p_2, \dots, p_{k-1}$, with \ $\lambda_2 : = \lambda_1 + p_1 -1$. So thanks to the induction hypothesis, there exists an open neighbourhood\ $U$ \  of \ $x_0$ \ in \ $X$ \ and an holomorphic map \ $g : U \to \mathcal{S}(\lambda_2, p_2, \dots, p_{k-1})$ such that \ $g^*(\mathbb{S}) $ \ is isomorphic to \ $( \mathbb{E}\big/\mathbb{F}_1)_{\vert U}$, where \ $\mathbb{S}$ \ is the canonical family parametrized by \ $\mathcal{S}(\lambda_2, p_2, \dots, p_{k-1})$.\\
 As the result is local near \ $x_0$, we may assume, up to thrink \ $U$ \ around \ $x_0$,  that \ $\mathbb{E}$ \ is given by an holomorphic \ $k-$thematic map \ $\varphi : U \to \Xi_{\lambda}^{(k-1)}$, which satisfies
 $$ \varphi(x) = s^{\lambda_k-1}.(Log\,s)^{k-1} + \psi(x)  $$
 where \ $\psi : U \to  \Xi_{\lambda}^{(k-2)}$ \ is holomorphic and \ $(k-1)-$thematic. The holomorphic map \ $g$ \ gives  in fact holomorphic maps \ $S_2, \dots, S_{k-1} : U \to 1 +b.\C[[b]]$ \ such that, if we define
 $$P_1 : = (a -\lambda_2.b).S_2^{-1}.(a - \lambda_3.b).S_3^{-1} \dots S_{k-1}^{-1}.(a - \lambda_k.b), $$
 \ $\A.P_1$ \ will be the annihilator of the standard generator \ $e$ \ of the canonical family parametrized by \ $\mathcal{S}(\lambda_2, p_2, \dots, p_k)$. So the generator  \ $g^*(e)$ \ of \ $\mathbb{E}\big/\mathbb{F}_1$ \ will satisfy also \ $P_1.g^*(e) = 0$. If we identify \ $\mathbb{E}_{\vert U}$ \ to \ $\mathbb{E}_{\varphi} \subset \Xi_{\lambda}^{(k-1)}$, we identify \ $\mathbb{E}\big/\mathbb{F}_1$ \ to a sub-sheaf of the quotient  sheaf  
 $$ \Xi_{X, \lambda}^{(k-1)}\Big/ \Xi_{X, \lambda}^{(0)}\  \simeq  \  \Xi_{X, \lambda}^{(k-2)}. $$
 So we may find \ $T_0, \dots, T_{k-1}$ \ local sections of \ $\mathcal{O}_X[[b]]$ \ such that the image of the section
 $$ \gamma : = \sum_{j=0}^{k-1} \ T_j.a^j.\varphi $$
 in \ $\mathbb{E}\big/\mathbb{F}_1$ \  co{\"i}ncides with \ $g^*(e)$. As \ $g^*(e)$ \ generates \ $\mathbb{E}\big/\mathbb{F}_1$ \ the section \ $T_0$ \ of \ $\mathcal{O}_X[[b]]$ \ must be  invertible near \ $x_0$. So \ $\gamma$ \ generates \ $\mathbb{E}$ \ near \ $x_0$ \ and satifies \ $P_1.\gamma \in \mathbb{F}_1$. As \ $\mathbb{F}_1 = \mathcal{O}_x[[b]].s^{\lambda_1-1}$ \ we may write
 $$ P_1.\gamma = \Theta.s^{\lambda_1-1} $$
 where \ $\Theta$ \ is a section of \ $ \mathcal{O}_X[[b]]$ \ in an open neighbourhood of \ $x_0$. The decomposition \ $E_{\lambda_1} = P_1.E_{\lambda_1} \oplus V_1$ \ (see [B.13-b] proposition 3.2.4) allows to write
 $$\Theta.s^{\lambda_1-1} = (P_1.\beta + S_1).s^{\lambda_1} $$
 where \ $\beta$ \ and \ $S_1$ \ are local sections respectiveley of \ $\mathcal{O}_X[[b]]$ \ and \ $\mathcal{O}_X\otimes V_1$. Moreover, as \ $\Theta$ \ is invertible, which means invertibility in \ $\mathcal{O}_X$ \ of its constant term in \ $b$, so is \ $S_1$. Then, up to multiplication by an holomorphic invertible function \ $I$ \   on an open neighbourhood of \ $x_0$, we may assume that the constant term of \ $S_1$ \ is \ $1$. Then \ $\tau : = I.(\gamma - \beta.s^{\lambda_1-1})$ \ is still a generator of \ $\mathbb{E}$ and satisfies
 $$ P_1\tau = S_1.s^{\lambda_1-1} \quad {\rm with} \quad S_1 \in \mathcal{O}_X \otimes V_1, S_1(0) = 1 .$$
 This gives \ $(a - \lambda_1.b)S_1^{-1}.P_1.\tau = 0 $. So we see that \ $\mathbb{E}$ \ is isomorphic to the pull back of the canonical family by the holomorphic map \ $f  : U(x_0) \to \mathcal{S}(\lambda_1, p_1, \dots, p_{k-1})$ \ given in the open neighbourhood \ $U(x_0)$ \ of \ $x_0$ \ by \ $S_1, \dots, S_k$. This isomorphism sends the local generator \ $\tau$ \ of \ $E$ \  on \ $ f^*(e)$ \ where \ $e$ \ is the standard generator of the canonical family parametrized by \ $ \mathcal{S}(\lambda_1, p_1, \dots, p_{k-1})$. So the local versality is proved.\\
 The proof of the global universality of this family under the assumption that any theme with these fundamental invariants has an unique canonical form is obvious. The fact that this is true when we assume \ $p_{j}\geq k-1 \quad \forall j \in [1,k-1]$ \ is proved in [B.13-b] theorem 4.3.1. $\hfill \blacksquare$
 
 \parag{Remark} If for some \ $t_{0}\in \mathcal{S}(\lambda_{1}, p_{1}, \dots, p_{k-1})$ \ there exists an open neighbourhood \ $U$ \ of \ $t_{0}$ \ such that any \ $t \in U$ \ has an unique canonical form, then the restriction to \ $U$ \ of the canonical family is locally universal at each point in \ $U$.

 \subsection{Examples.}
 
We shall give a examples of  \ $[\lambda]-$primitive themes of rank \ $3$ \ for which there does not exists, even locally, an universal family.\\

Let  the fundamental invariants be \ $\lambda_1 : = \lambda, p_1 = p_2 = 1$;  so we have  
$$\lambda_{1}  = \lambda_{2} = \lambda_{3} = \lambda, \quad q_1 = p_1 + p_2 =  2 \quad {\rm and} \quad q_2 = p_2 = 1.$$ 
Remark that the corresponding themes are special (i.e. $\lambda_{1} = \lambda_{2} = \cdots = \lambda_{k}$).
The corresponding parameter space for the canonical family is \ $\mathcal{S}(\lambda,1,1) \simeq (\C^{*})^{2}\times \C$ \ and we shall denote 
 $$E(\alpha,\beta,\gamma) = \A\big/\A.(a - \lambda.b).(1 + \beta.b + \gamma.b^2)^{-1}. (a - \lambda.b).(1 + \alpha.b)^{-1}.(a - \lambda.b) $$
 the theme associated to the parameter \ $(\alpha,\beta,\gamma) \in  (\C^{*})^{2}\times \C .$

\begin{prop}\label{ex. 1}
\begin{enumerate}
\item For \ $\alpha \not= \beta$ \ and any \ $\gamma$ \ the theme \ $E(\alpha,\beta,\gamma)$ \ is not invariant, and its canonical form is not unique. 
\item For \ $\alpha = \beta \not= 0$ \ the theme \ $E(\alpha,\alpha,\gamma)$ \ is invariant. But there does not exist an universal family near each such invariant  theme.
\item For any \ $\alpha \not= \beta$ \ and any \ $\gamma$, there exists a locally universal family with dimension \ $2$ \ smooth parameter space.
\end{enumerate}
\end{prop}

The proof will use  the following two lemmas.

\begin{lemma}\label{isom.}
Let  \ $\alpha,\beta, \gamma \in \C, \alpha.\beta \not= 0 $; then  the rank \ $3$ \ theme \ $E(\alpha,\beta,\gamma)$ \ is defined as follows :
\begin{align*}
&  (a - \lambda.b).e_3 = (1 + \alpha.b).e_2 \\
&  (a - \lambda.b).e_2 = (1 + \beta.b + \gamma.b^2).e_1 \\
& (a - \lambda.b).e_1 = 0 .
\end{align*}
For \ $\beta \not= \alpha$, \ $E_{\alpha,\beta,\gamma}$ \ is isomorphic to \ $E_{\alpha,\beta, 0}$ \ for each \ $\gamma \in \C$.\\
For  \ $\beta = \alpha$, the themes \ $E_{\alpha,\alpha,\gamma}$ \ and \ $E_{\alpha,\alpha,\gamma'}$ \ are isomorphic if and only if \ $\gamma = \gamma'$.
\end{lemma}

Note that this lemma already implies the non uniqueness of the canonical form of \ $E(\alpha,\beta,\gamma)$ \ for \ $\alpha \not= \beta$, and also that near any such theme, the canonical family is not locally universal. 

\parag{Proof} We look for a \ $\C[[b]]-$basis \ $\varepsilon_3, \varepsilon_2, \varepsilon_1$ \ of \ $E_{\alpha,\beta,\gamma}$ \ satisfying the following conditions :
\begin{align*}
& \varepsilon_3 = e_3 + U.e_2 + V.e_1, \quad {\rm with} \quad U,V \in \C[[b]] \tag{0}\\
& (a - \lambda.b).\varepsilon_3 = (1 + \alpha.b).\varepsilon_2 \tag{1} \\
& (a - \lambda.b).\varepsilon_2 = (1 + \beta.b + \gamma'.b^2).\varepsilon_1 \tag{2}\\
& (a - \lambda.b).\varepsilon_1 = 0. \tag{3}
\end{align*}
We know that \ $\alpha$ \ and \ $\beta$ \ are determined by the isomorphism class of the theme \ $E(\alpha,\beta,\gamma)$:  as we have \ $p_1 = p_2 = 1$ \ they are the parameters of the rank \ $2$ \ themes \ $E\big/F_1$ \ and  \ $F_2$ \ respectively .\\
Remark also that \ $(3)$ \ implies that \ $\varepsilon_1 = \rho.e_1$ \ with \ $\rho \in \C^*$.\\
Let us compute the conditions for \ $U$ \ and \ $V$ :
\begin{align*}
& (a - \lambda.b).\varepsilon_3 =   (1 + \alpha.b).e_2 + b^2.U'.e_2 + U. (1 + \beta.b + \gamma.b^2).e_1 + b^2.V'.e_1 \\
&\qquad\qquad   = (1+\alpha.b).\varepsilon_2
\end{align*}
and so, we have \ $ \varepsilon_2 = Z.e_2 + T.e_1$ \ with \ $ Z = (1 +\alpha.b)^{-1}.(1 + \alpha.b + b^2.U') $ , and 
\begin{equation*}
 (1 + \alpha.b).T = U.(1 + \beta.b + \gamma.b^2) + b^2.V'. \tag{4}
\end{equation*}
Then we get
\begin{align*}
& (a - \lambda.b).\varepsilon_2 = Z.(1 + \beta.b + \gamma.b^2).e_1 + b^2.Z'.e_2 + b^2.T'.e_1 \\
& \qquad \qquad = (1 + \beta.b + \gamma'.b^2).\rho.e_1
\end{align*}
and this implies \ $Z' = 0 $ \ and as  \ $Z = 1 + (1 + \beta.b)^{-1}.b^2.U'$ \ we must have \ $U \in \mathbb{C}$, and \ $Z = 1$.  The relation \ $(2)$ \ implies now, as  \ $\varepsilon_2 = e_2 + T.e_1$
$$  (1 + \beta.b + \gamma.b^2).e_1 + b^2.T'.e_1 =  (1 + \beta.b + \gamma'.b^2).\rho.e_1 .$$
Then we have \ $\rho = 1$ \ and \ $ T' = \gamma' - \gamma $. \\
Looking at the constant terms in \ $(4)$ \ we obtain  \ $T = U + (\gamma' - \gamma).b$.\\
But \ $(4)$ \ implies also
\begin{equation*}
 \alpha.U + \gamma' - \gamma = U.\beta \qquad \qquad {\rm and}  \quad
 U.\gamma + V' = \alpha.(\gamma - \gamma') \tag{5} 
\end{equation*}
So, for \ $\alpha \not= \beta $ \ we will have
\begin{equation*}
U = \frac{\gamma - \gamma'}{\alpha - \beta} \qquad {\rm and} \qquad  V = V_0 + \frac{\gamma'- \gamma}{\beta - \alpha}.\big(\alpha.(\beta-\alpha) - \gamma\big).b . \tag{6}
\end{equation*}
If \ $\beta = \alpha$, the relation \ $(5)$ \ gives  \ $\gamma = \gamma'$. $\hfill \blacksquare$

\bigskip

The next lemma shows that for \ $\alpha = \beta \not= 0$ \  $E(\alpha, \beta,\gamma)$ \  is an invariant theme for any \ $\gamma \in \C$.

\begin{lemma}\label{Inv. alpha = beta}
For \ $\alpha \not= 0$ \ the (a,b)-module \  $E_{\alpha,\alpha,\gamma}$ \ is an invariant rank \ $3$ \ theme.
\end{lemma}

\parag{Proof} It is enough to find  \ $x : = e_2 + U.e_1$  \ such that
$$(a-\lambda.b)(1+\alpha.b)^{-1}(a -\lambda.b).x = 0 , $$
where \ $U \in \C[[b]]$. \\
As \ $F_2$ \ is a theme, the kernel of \ $(a - \lambda.b)$ \ in \ $F_2$ \ is \ $\C.e_1$. So \ $x$ \ must satisfies
$$ (a - \lambda.b)x = \rho.(1 + \alpha.b).e_1 $$
and this implies that \ $U$ \ is solution of the equation
$$ (1 + \alpha.b+ \gamma.b^2) + b^2.U' = \rho.(1 + \alpha.b) .$$
So we conclude that we must have \ $\rho = 1$ \ and  \ $U = -\gamma.b + cste$. So we obtain a solution \ $ x : = e_2 -\gamma.b.e_1$. $\hfill \blacksquare$\\

For  \ $ (\alpha,\beta) \in (\C^*)^2, \alpha \not= \beta $ \ denote \ $E(\alpha, \beta) : = E(\alpha,\beta,0)$ \ the rank \ $3$ \ theme defined by 
 $$E(\alpha,\beta) : = \A\big/\A.(a -\lambda.b)(1+\beta.b)^{-1}(a - \lambda.b)(1 + \alpha.b)^{-1}(a - \lambda.b).$$

\begin{lemma}\label{Univers.}
The family  \ $E(\alpha,\beta)_{(\alpha,\beta,0)\in X}$ \ is universal near any point in \\
 $X : = (\C^*)^2 \setminus \{\alpha = \beta\} $.
\end{lemma}

\parag{Proof} Thanks to the lemma \ref{isom.} it is enough to prove that for \ $\alpha \not= \beta$ \  the restriction of the canonical family to \ $X \subset \{\gamma = 0 \}$ \ is locally versal. Let \ $\mathbb{E}$ \ be the sheaf on \ $X$ \ of \ $\mathcal{O}_X-(a,b)-$modules associated to the holomorphic family \ $E(\alpha,\beta,0)$. It is enough to prove that the holomorphic map
$$\pi :  X \times \C \to X $$
defined by \ $\pi(\alpha,\beta, \gamma) = (\alpha,\beta)$ \ is such that the sheaf \ $\pi^*(\mathbb{E})$ \ is isomorphic as a sheaf of \ $\mathcal{O}_X-(a,b)-$modules, to the sheaf associated to the canonical family restricted to \ $X \times \C$. As the map is given by the versality of the canonical family, so  it is enough to construct its inverse. The inverse  we are looking for this map is given by the computation in the lemma \ref{isom.}, because  via the formulas \ $(6)$ \ in the case \ $\gamma' = 0$ \  we may send the generator \ $e_3$ \ of \ $E(\alpha,\beta,\gamma)$ \  to  \ $\varepsilon_3 : = e_3 + U.e_2 + V.e_1$, which is the generator of the \ $E(\alpha,\beta)$ \ for \ $\alpha \not= \beta$. $\hfill \blacksquare$

\parag{Proof of the proposition \ref{ex. 1}} The fact that the invariant themes in this canonical family are exactely the \ $E_{\alpha,\alpha,\gamma}$ \ is proved in lemmas \ref{isom.} and \ref{Inv. alpha = beta}.\\
The canonical family  \ $(E_{\alpha,\beta,\gamma})_{(\alpha,\beta,\gamma) \in S(\lambda_1,  p_1 = p_2 =1)}$ \ is holomorphic and versal near any point thanks to the theorem \ref{versal univ. 2}. Let assume that we have find an holomorphic family  \ $(E_y)_{y \in Y}$ \  which is locally universal around a theme \ $ E(\alpha_0,\alpha_0,\gamma_0) \simeq E_{y_0}$, where \ $Y$ \ is a reduced complex space that we may assume to be embedded in \ $\C^N$ \ near \ $y_0$. Let  $\varphi : \Omega \to Y \hookrightarrow \C^N$ \ the classifying map for the canonical family on an open set  \ $\Omega$ \ de \ $(\alpha_0,\alpha_0,\gamma_0) \in (\C^*)^2\times \C $. As for \ $\alpha \not= \beta$ \ the isomorphy class of the theme $E_{\alpha,\beta,\gamma}$ \  does not depend on \ $\gamma$, thanks to the lemma \ref{isom.}, we shall have \ $\frac{\partial \varphi}{\partial \gamma} \equiv 0$ \ on the open set \ $\{\alpha \not= \beta\}$ \ of \ $\Omega$. This implies that \ $\varphi$ \  does not depends on \ $\gamma$ \ for all  \ $\alpha$ \ near enough to \ $\alpha_0$, and for all \ $\gamma,\gamma'$ \ near enough from \ $\gamma_0$. This contradicts the lemma \ref{isom.}. \\
The last statement of the proposition is given by lemma \ref{Univers.}. $\hfill \blacksquare$

 \parag{Another example without universal family}
 
 We give here an example of {\bf non special} rank \ $4$ \ $[\lambda]-$primitive themes without universal family.
 Fix the fundamental invariants to be \ $\lambda_1 > 3$ \ and \ $p_1 = 3, p_2 = p_3 = 2$. So the canonical family is defined by the following relations :
 
 \begin{align*}
 & (a - (\lambda_1+4).b).e_4 = S_3.e_3 \quad S_3 : = 1 + \alpha.b^2  \quad\quad (q_3 = p_3 = 2) \\
 & (a - (\lambda_1+3).b).e_3 = S_2.e_2 \quad  S_2 : = 1 + \beta.b + \gamma.b^2 \quad\quad (q_2 = p_2 = 2) \\
 & (a - (\lambda_1+2).b).e_2 = S_1.e_1 \quad S_1 : = 1 + \delta.b + \varepsilon.b^2 + \theta.b^3 \quad\quad  (q_1 = p_1 = 3)\\
 & (a - \lambda_1.b).e_1 = 0.
 \end{align*}
 with the condition  \ $ \alpha.\gamma.\theta \not= 0 $.\\
 
 We look for the invariant themes in this family. So we look for an element \ $x \in F_3 \setminus F_2$ \ such that \ $P_0.x = 0$. The existence of such an \ $x$ \  is equivalent to the existence of a rank \ $3$ \ endomorphism.\\
 Using proposition 3.2.2 of [B.13-b], the condition \ $P_0.x = 0$ \ is in fact equivalent to the condition \ $P_1.x = 0$ \ because we ask \ $x$ \ to be in \ $F_3$. This gives the equality \ $(a - \lambda_2.b).S_2^{-1}.P_2.x = 0$.
 The kernel in \ $E$ \ of \ $(a - \lambda_1.b)$ \ is \ $\C.e_1$, so we obtain the equation
 \begin{equation*}
 P_2.x = \rho.S_2.b^2.e_1 \tag{1} 
 \end{equation*}
 It gives \ $(a - \lambda_3.b).S_3^{-1}.P_3.x = \rho.S_2.b^2.e_1$ \ which implies that the class of \  \ $S_3^{-1}.P_3.x$ \ in \ $E\big/F_1$ \  is in the kernel of \ $(a - \lambda_3.b)$ \ which is  \ $\C.b^{\lambda_3-\lambda_2}.e_2 \quad modulo \  F_1$. So we may put :
 \begin{equation*}
 P_3.x = \sigma.S_3.b.e_2 + S_3.T.e_1 \tag{2}
 \end{equation*}
and the equation \ $(1)$ \ implies 
 \begin{align*}
 &   P_2.x = (a - \lambda_3.b)( \sigma.b.e_2 + T.e_1) = \rho.S_2.b^2.e_1 \\
 & \qquad = \sigma.b.S_1.e_1 + (b^2.T' - 3b.T).e_{1} =  \rho.S_2.b^2.e_1. 
 \end{align*}
So we get \ $ b^2.T' - 3.b.T = \rho.S_2.b^2 - \sigma.S_1.b $ \ and after simplification by \ $b$ \ this gives 
 \begin{equation*}
 b.T' - 3.T = \rho.S_2.b - \sigma.S_1. \tag{3}
 \end{equation*}
 Then we conclude, as the coefficient of \ $b^{3}$ \ must be zero in the right handside, that we must have \ $\rho.\gamma = \sigma.\theta$. As \ $\gamma.\theta \not= 0$, the condition \ $\rho \not= 0 $ \ is equivalent to \ $\sigma \not= 0$.\\
 
 We shall have a solution \ $T$ \ by choosing \ $\sigma : = \rho.\gamma/\theta$, and it satisfies   \ $-3T(0) = -\sigma$, that is to say \ $T(0) = \sigma/3 =  \rho.\gamma/3\theta$.\\
 Define now  \ $x : = X.e_3 + U.e_2 + V.e_1 $. The relation \ $(2)$ \ gives 
 $$ (a - \lambda_4.b).x = (b^2.X' - b.X).e_3 \quad {\rm modulo} \ F_2 $$
 and so  \ $X = \tau.b$ \ with \ $\tau \in \C^*$ \ as we assume \ $x \not\in F_2$. Then we get
 \begin{align*}
& (a - \lambda_4.b).x = (\tau.b.S_2 + b^2.U' - 2b.U).e_2 + (b^2.V' - 4b.V + U.S_1).e_1 \\
& \qquad    = \sigma.S_3.b.e_2 + S_3.T.e_1 
\end{align*}
and so the equations
\begin{align*}
& b.U' - 2U = \sigma.S_3 - \tau.S_2  \tag{5} \\
& b^2.V' - 4b.V = S_3.T - U.S_1 \tag{6}
\end{align*}
The first one implies that \ $\tau.\gamma = \sigma.\alpha$ \ and forces \ $U(0) = (\tau - \sigma)/2= \sigma.(\alpha/\gamma -1)/2  $; the second one imposes
$$ T(0) = U(0) .$$
Then we obtain
$$ T(0) = \sigma/3  = \sigma(\alpha/\gamma - 1)/2 . $$
So, if \ $ \alpha/\gamma \not= 1 + 2/3 = 5/3 $ \ we must have  \ $\sigma = 0$ \ and this imposes \ $\tau = 0$ \ and this is incompatible with our assumption \ $x \in F_3 \setminus F_2$.\\
Then for  \ $3\alpha \not= 5\gamma $ \ the theme is not invariant, and as  \ $E/F_1$ \ is invariant thanks to the theorem 4.3.1 of [B.13-b] (here the rank is \ $3$ \ and each \ $p_i$ \ is at least \ $2$), the canonical form will not be unique thanks to lemma 4.3.4 of [B.13-b].\\
Now for \ $3\alpha = 5\gamma $ \ we find a solution \ $U$ \ of the equation \ $(5)$ \ and then a solution \ $V$ \ for \ $(6)$ \ when we choose a solution \ $T$ \ of \ $(3)$. So along the hyperplane  \ $3\alpha = 5\gamma $ \ the themes are invariant.  The situation is similar to the previous example.\\

Let us show, to conclude the proof, that when \ $3\alpha = 5\gamma $ \ the coefficient  \ $\varepsilon$ \ is not relevant to determine the isomorphism class of \ $E$. As we know that \ $\theta$ \ is determined by the isomorphy class of \ $E$, it is enough to prove that \ $\delta$ \ is also determined by the isomorphy class of \ $E$.\\
 Because, up to an homothetie, any automorphism of \ $E$ \ is obtained  by sending \ $e_4$ \ to \ $e_4 + y$ \ where \ $y \in F_3$  \ and satisfies \ $P_1.y \in F_1$, to prove that \ $\delta$ \ is determined by the isomorphy class of \ $E$, it is enough to prove that any such \ $y$ \ satisfies \ $P_1.y \in b^2.F_1$.\\
But  \ $P_1.y \in F_1$ \  shows that if we define \ $\varphi(e) = y$ \ this produces an \ $\A-$linear map \ $E\big/F_1 \to F_3\big/F_1$. As the rank is \ $\leq 2$, $\varphi(F_2\big/F_1) = 0$ \ and so  \ $P_2.e \in F_2$ \ implies \ $P_2.y \in F_1$. Then we have 
$$ (a - \lambda_3.b).S_3^{-1}.P_2.y \in F_1 .$$
The kernel of \ $(a - \lambda_3.b)$ \ acting on  \ $F_3\big/F_1$ \ is \ $\C.b.e_2 + F_1$ \ as \ $\lambda_3 - \lambda_2 = 1$. Then we may write
$$ P_3.y = \rho b.S_3.e_2 + S_3.T.e_1 $$
which gives 
$$ P_2.y = \big[ \rho b.S_1 + b^2.T' - 3b.T \big].e_1 \in b.F_1 .$$
So  \ $P_1.y \in b^2.F_1$ \ is proved. \\
It is not necessary to prove that \ $\varepsilon$ \ can move by some isomorphism because  this is consequence of the lemma 4.3.4 in [B.13-b]. $\blacksquare$

\parag{Conclusion} So in this example we have the following properties, analogous of the properties in the previous rank \ $3$  example:
\begin{enumerate}
\item For \ $3\alpha \not= 5\gamma$ \ the themes are not invariant.
\item For \ $3\alpha = 5\gamma$ \ each theme is invariant but there is no local universal family around such an isomorphy class.
\item For each theme \ $E$ \  such that \ $3\alpha \not= 5\gamma$ \ there exists a local universal family obtained by the restriction of the canonical family to an open neighboorghood of \ $E$ \ in the hyperplane \ $\{\varepsilon = \varepsilon_0\}$.\\
\end{enumerate}

\section{Bibliography}

\begin{itemize}

\item{[A-G-V]} Arnold, V. \  Goussein-Zad\'e, S. \ Varchenko,A. {\it Singularit\'es des applications diff\'erentiables}, volume 2,  \'edition MIR, Moscou 1985.

\item{[B.93]} Barlet, D. {\it Theory of (a,b)-modules I}, in Complex Analysis and \\ Geometry, Plenum Press New York (1993), p.1-43.

\item{[B.95]} Barlet, D. {\it Theorie des (a,b)-modules II. Extensions}, in Complex \\ Analysis and Geometry, Pitman Research Notes in Mathematics, Series 399 Longman (1997), p.19-59.

\item{[B. 05]}  Barlet, D. \textit{Module de Brieskorn et formes hermitiennes pour une singularit\'e isol\'ee d'hypersurface}, Revue Inst. E. Cartan (Nancy) vol. 18 (2005), p.19-46.  

\item{[B.08]} Barlet, D.  \textit{Sur certaines singularit\'es d'hypersurfaces II},   Journal of \\ Algebraic 
  Geometry 17 (2008), p.199-254. 

\item{[B.09]} Barlet, D. {\it P\'eriodes \'evanescentes et (a,b)-modules monog\`enes},  Bulletino U.M.I. (9) II (2009), p.651-697.


\item{[B.13-b]} Barlet, D. {\it The theme of a vanishing period} to appear.

 \item{[Br.70]} Brieskorn, E. {\it Die Monodromie der Isolierten Singularit{\"a}ten von \\ Hyperfl{\"a}chen}, Manuscripta Math. 2 (1970), p.103-161.
 
 \item{[H.99]} Hertling, C. {\it Classifying spaces for polarized mixed Hodge structures and for Brieskorn lattices}, Compositio Math. 116 (1999), p.1-37.

\item{[K.09]} Karwasz, P. {\it Self dual (a,b)-modules and hermitian forms}, Th\`ese de doctorat de l'Universit\'e H. Poincar\'e (Nancy 1) soutenue le 10/12/09, p.1-57.

\item{[K.76]} Kashiwara, M.  {\it b- function and holonomic systems}, Inv. Math. 38 (1976) p.33-53.

\item{[M.74]} Malgrange, B.   {\it Int\'egrales asymptotiques et monodromie}, Ann. Sci. Ec. Norm. Sup. 4 (1974) p.405-430.

  \item{[M.75]} Malgrange, B. {\it Le polyn\^ome de Bernstein d'une singularit\'e isol\'ee}, in Lect. Notes in Math. 459, Springer (1975), p.98-119.

\item{[S.89]} Saito, M. {\it On the structure of Brieskorn lattices}, Ann. Inst. Fourier 39 (1989), p.27-72.

\item{[S.91]} Saito, M. {\it Period mapping via Brieskorn modules}, Bull. Soc. Math. France 119 (1991), p.141-171.

\end{itemize}

\end{document}